\documentclass{llncs}
\usepackage{llncsdoc}

\begin{document}

\title{System Description: Russell - A Logical Framework for Deductive Systems}
\author{Dmitry Vlasov\inst{1}}
\institute{Sobolev Institute of Mathematics, SB RAS \email{vlasov@math.nsc.ru}
\footnote{The research is partially supported by Russian Foundation for Basic Researches (Grant No. 17-01-00531)}}

\maketitle

\begin{abstract}
Russell is a logical framework for the specification and implementation of 
deductive systems. It is a high-level language with respect to Metamath language \cite{Metamath}, so
inherently it uses a Metamath foundations, i.e. it doesn't rely on any particular
formal calculus, but rather is a pure logical framework. The main difference with Metamath
is in the proof language and approach to syntax: the proofs have a declarative form, i.e.
consist of actual expressions, which are used in proofs, while syntactic grammar rules are
separated from the meaningful rules of inference.

Russell is implemented in c++14 and is distributed under GPL v3 license.
The repository contains translators from Metamath to Russell and back. Original Metamath theorem base 
(almost 30 000 theorems) can be translated to Russell, verified, translated back to
Metamath and verified with the original Metamath verifier. Russell can be downloaded from  the repository \verb|https://github.com/dmitry-vlasov/russell|
\end{abstract}

\section{Introduction}
Recently the ambitious QED project \cite{QED} has celebrated its 20 year anniversary, while 
the claimed goals of this project are still far from being reached. Several papers \cite{17_provers}, \cite{QED_revisited}, 
\cite{QED_Hammering}, \cite{QED_reloaded}, \cite{Challenge}
addressing the history of QED clearly state that yet there is no computer language, which has all the expected features of a QED system. Summarizing these papers, it can be said that the 
major barriers of QED are:
\begin{itemize}
 \item the 'Balcanization' of QED-like systems, i.e. when there are different languages with different foundations and there is no simple way to share formalized proofs between them \cite{QED_reloaded}
 \item the lack of a powerful automation, which could seriously reduce end user efforts to prove a theorem \cite{QED_Hammering}
 \item the difference between the standard practical mathematical language, which is used in papers and textbooks, and the formal language of QED-like systems \cite{QED_revisited}
\end{itemize}
Russell is another system, which is intended to reach QED goals. It was developed to address all these obstacles, and some of them are to some extent removed in the Russell design.

\section{The Russell System}

The Russell system is a general purpose framework for definition and 
usage of different formal deductive systems. The variety of formal systems,
which can be represented in Russell, is quite wide, although limited.
For example, non-monotonic deductive systems cannot be given in Russell.
As a high-level language wrt Metamath \cite{Metamath}, Russell inherits
its foundations, which are close to the notion of Post's canonical system \cite{Post}.
The approach of Metamath to syntax of expressions is more general than that of
Russell: syntactic rules in Metamath are indistinguishable from the
meaningful rules of inference and can have meaningful (essential, in Metamath
terms) premises. In Russell, the grammar of expressions must be context free by 
the language definition. In all other aspects Metamath and Russell share the 
same foundations.

\subsection{Pure Logical Framework}

What distinguishes Metamath and Russell from other logical frameworks is purity:
their deduction engines don't use any built-in logic (in the form of axioms and inference rules). 
The deduction used in Metamath and Russell is concerned with making a proper substitution, 
applying it to a certain expression and checking the coincidence of the result with some other
expression. From the very general point of view, such 
logic-neutrality is a good property, because we want to avoid the situation when 
some formal deductive system has an a-priory better fitness to the language
than the other 
because of the propinquity with the underlying logic of a logical framework. 
Similar arguments are mentioned in the paper \cite{QED_reloaded}, where the statement 
of \textit{foundational pluralism} is given. In fact, the property of
logic-neutrality (or pureness) is crucial to the desired foundational pluralism, 
because it gives a uniform core language for a vast variety of deductive systems. 
Thus, 'Balkanization' of formal mathematics could be at least based on the same
language (although the sharing of proofs between different foundations is still challenging).

What distinguishes Russell from Metamath:
\begin{itemize}
 \item Expression grammar syntax rules are separated from the general logical 
 inference rules, so the grammar is guaranteed to be context free.
 \item The Russell proof language is more human-friendly than the Metamath proof 
 language.
 \item Russell uses a special syntactic construct for definitions, while in 
 Metamath definitions are simply axioms with labels, starting with 'df' prefix.
\end{itemize}

\subsection{Syntax of expressions in Metamath and Russell}

The syntax of expressions in Metamath is not distinguished from the
syntax of axioms and inference rules. The latter assertions are treated differently inside proofs, but on the syntax level there is no difference between these two kinds of assertions.

For example, 
the definition of well-formed-formula with $\rightarrow$ and $\neg$ logical
connectives in Metamath looks like:
\begin{verbatim}
${
    wph $f wff ph $.
    wi $a wff ( ph -> ps ) $.
$}
${
    wph $f wff ph $.
    wps $f wff ps $.
    wi $a wff ( ph -> ps ) $.
$}
\end{verbatim}
and the definition of an axiom of the Hilbert-style propositional calculus looks like:
\begin{verbatim}
${
    wph $f wff ph $.
    wps $f wff ps $.
    ax-1 $a |- ( ph -> ( ps -> ph ) ) $.
$}
\end{verbatim}
Technically, one could add substantial (i.e. essential, in the Metamath terminology) hypotheses
to the syntactic grammar rules, but of course, it would make little sense, since the
common practice is to use context-free grammars for a language.

The same syntactic rules and axioms in Russell look like:
\begin{verbatim}
rule wn (ph : wff) {
    term : wff = # ¬ ph ;;
}
rule wi (ph : wff, ps : wff) {
    term : wff = # ( ph → ps ) ;;
}
axiom ax-1 (ph : wff, ps : wff)  {
    prop 1 : wff = |- ( ph → ( ps → ph ) ) ;;
}
\end{verbatim}
Context-freeness here is obligatory by design.

\subsection{Proof Language Syntax in Metamath and Russell}

Proofs in Metamath are simply sequences of labels, which mark up assertions. 
Some of them mark up syntax-forming assertions, which are actually grammar rules, 
other are used to mark meaningful axiomatic or provable assertions. These sequences
form a reversed polish notation (RPN) for a program, which evaluates the final 
expression. The proof verification procedure consists of two steps: computing
the expression from a proof and comparing this expression with the proposition 
of a proved assertion.

Thus, the only way to understand the Metamath proof is to compute it as an RPN
program, which demands a proof assistant (except for some trivial cases). 
An example of the Metamath provable assertion:
\begin{verbatim}
${
    syl.1 $e |- ( ph -> ps ) $.
    syl.2 $e |- ( ps -> ch ) $.
    syl $p |- ( ph -> ch ) $= .... $.
$}
${
    $v ph ps ch $.
    wph $f wff ph $.
    wps $f wff ps $.
    wch $f wff ch $.
    a1d.1 $e |- ( ph -> ps ) $.
    a1d $p |- ( ph -> ( ch -> ps ) ) $=
      wph wps wch wps wi a1d.1 wps wch ax-1 syl $.
$} 
\end{verbatim}
Here the sequence 
\begin{verbatim}
wph wps wch wps wi a1d.1 wps wch ax-1 syl 
\end{verbatim}
is an RPN program (Metamath proof) which produces the statement 
\begin{verbatim}
( ph -> ( ch -> ps ) ) 
\end{verbatim}
from the (essential, in Metamath terminology) hypothesis \verb|( ph -> ps )|. This short
proof can be understood by a human directly, but more complex proofs demand a proof 
assistant, like a metamath program, which would generate the stack trace of 
execution of the proof as an RPN program. 

The Russell proof language instead uses a much more user-friendly declarative
form, in which a proof is a sequence of steps, each of which is an expression 
accompanied with references to the appropriate assertions and premises, which are
previously proved expressions or expressions from the premise list of a 
currently proved theorem. Technically, this sequence of steps is a stack
trace obtained from executing a proof in an RPN form, which is stripped 
off all intermediate steps related to the syntax formation. When a Russell proof
is verified or translated to Metamath, these intermediate steps are restored from
the syntax tree of an expression and unificators. 


Example of the Russell proof for the same assertion:
\begin{verbatim}
theorem syl (ph : wff, ps : wff, ch : wff)  {
  hyp 1 : wff = |- ( ph -> ps ) ;;
  hyp 2 : wff = |- ( ps -> ch ) ;;
  -----------------------
  prop 1 : wff = |- ( ph -> ch ) ;;
}
proof of syl { ... }
\end{verbatim}
\begin{verbatim}
theorem a1d (ph : wff, ps : wff, ch : wff)  {
  hyp 1 : wff = |- ( ph -> ps ) ;;
  -----------------------
  prop 1 : wff = |- ( ph -> ( ch -> ps ) ) ;;
}
proof of a1d {
  step 1: wff = axm ax-1 () |- ( ps -> ( ch -> ps ) ) ;;
  step 2: wff = thm syl (hyp 1, step 1) |- ( ph -> ( ch -> ps ) ) ;;
  qed prop 1 = step 2 ;
}
\end{verbatim}

Thus, the Russell proof language is natural and
simple, which imposes minimum restrictions and allows for user-defined grammars for expressions. Together with the possibility to use the conventional set theory it bridges the gap between the computer-based mathematics and the common practice mathematics from 
textbooks.

\subsection{Definitions in Metamath and Russell}

One of the most important features of any computer-based deductive system (framework) 
is safety and reliability \cite{Proof_audit}: 

  \emph{to what extent one can trust computer proofs} ? 

\noindent
Reliability of a formal system is a complex subject, which involves several aspects. One of
these aspects is the size of the axiomatic base used by a theory. If it is 
large then there can be some (unintentionally) hidden inconsistency inside of it, and if so,
this will lead to the triviality of the whole theory (in case of an explosive logic). On the other hand,
if the set of the true axioms is small and well-known (like some variant of ZF set theory) then
its degree of reliability is very high. 

If each definition is introduced as a new axiom, like it is done in Metamath, then the number 
of axioms increases fast as the theory grows, and at some point there is no
guarantee that all of these axioms are consistent. To address this issue, definitions in Russell are
introduced as a special syntactic construct and certain properties are checked for each 
definition to ensure that adding the underlying axiom will give a conservative extension of the
theory. Conservativity here means that if something can be proved with the help of some
definition, it can also be proved without it. This property is strictly proven, so it gives
some more certainty about correctness of Russell theories.

Example of a definition in Russell:
\begin{verbatim}
definition df-or (ph : wff, ps : wff)  {
    defiendum : wff = # ( ph \/ ps ) ;;
    definiens : wff = # ( -. ph -> ps ) ;;
    -----------------------
    prop : wff = |- ( defiendum <-> definiens ) ;;
}
\end{verbatim}

The fact that all Russell sources can be translated back to Metamath and checked with its
original proof checker shows that Russell at least as reliable as Metamath. 
Moreover, the declarative format of proofs in Russell makes it possible for a human to do
an independent verification of proofs. Of course, it would look strange to exploit human ability of checking a formal proof in a QED system, but still, from the philosophical point of view,
human understanding is an ultimate judge, and is very important. 

\section{Implementation}
Currently the Russell language is implemented as a translator from / to Metamath and is written in c++14. The Russel repository includes test scripts, which run a chain of translations (here abbreviation MM stands for Metamath, SMM for simplified Metamath): 

\verb|MM -> SMM -> Russell -> SMM -> MM| 

The translation of the whole Metamath base (about 30 000 theorems) is rather fast in all directions. The most problematic from the performance point of view is expression parsing in Russell. Metamath uses an explicit construction of expressions in proofs, so it does not require
any parsing or unification algorithm when checking its source. Unlike Metamath, 
Russell must parse expressions in order to get syntax trees and such parsing takes the most of time in comparison to all other steps, like unification or translation.

One of the features implemented in the toolchain of translators is that it can automatically 
divide the original Metamath source (the single file of almost 150 megabytes) into reasonably
small parts following the internal layout inside the source file. After breaking it into pieces, 
one can browse the source file tree with the standard desktop navigation tools and watch source files in the standard desktop editors without the necessity to handle a single 150 Mb file.

Russell is not considered by the author as an experimental or model language. It is supposed to be a useful, universal, and convenient tool for all kinds of activity in the field of formal deduction. To achieve this, the language of implementation (c++), quality of source code, efficiency of algorithms, and usefulness to the end user are of a great importance. Russell implementation should be able to work with hundreds of thousands of assertions in a reasonable time, which is the subject of ongoing research and development.

\section{Conclusion and Future Work}
The Russell logic framework is a robust, fast, and reliable general purpose tool for representation of formal deductive systems. The Russell language is designed to be simple and easy-to-learn and provides proofs in a declarative form (which is a standard practice in informal mathematical texts). 

Essentially, powerful automation is the only one blocking property left in the list of the desirable QED features. Therefore, the next challenge is to implement the automated proving feature so that the process of
formal proof design would be easier. The first step to making an automated 
proof engine for Russell has been already taken and it showed the potential feasibility of such a goal. Since the proof search space suffers from an extreme 
combinatorial explosion, some extraordinary means are needed to cope with it. Standard techniques will not work here, since the nature of the underlying deductive system is apriori unknown to the prover (which is a consequence of the logic-neutrality property). For example, in general we can not assume, that the underlying logic is cut free (and actually it is not in the case of the Metamath theorem base).

To create a powerful prover for Russell we plan to use advanced machine earning techniques to make the prover use the experience of the already proven theorems. Ideally, it should be able to generate human-like proofs, formed as a combination of previously obtained proofs. 

The other important goal is to support importing of other bases of formalized mathematics into Russell. Some successful
attempts of importing HOL theorem base into Metamath have been already undertaken \cite{HOL_MM}, so there is a hope that it would
possible to join a large part of the already formalized mathematical knowledge under a common logical framework, but with different foundations. A more ambitious goal is to join these bases upon a common foundation.

\end{document}